\documentclass[11pt,3p,times]{elsarticle}

\usepackage{srcltx}
\usepackage{enumitem}
\usepackage{amssymb, color}
\usepackage{xcolor}
\usepackage{hyperref}
\usepackage{amsmath}
\usepackage{pifont}
\usepackage[utf8]{inputenc}

\makeatletter
\def\LaTeX{\leavevmode L\raise.42ex
    \hbox{\kern-.3em\size{\sf@size}{0pt}\selectfont A}\kern-.15em\TeX}
\makeatother

\newcommand{\BibTeX}{{\rm B\kern-.05em{\sc
          i\kern-.025emb}\kern-.08em\TeX}}
\makeatletter
\def\@currentlabel{2.1}\label{e:dispaa}
\def\@currentlabel{2.21}\label{e:dispau}
\def\@currentlabel{2.22}\label{e:dispav}
\def\@currentlabel{2.23}\label{e:dispaw}
\def\@currentlabel{2.24}\label{e:dispax}
\def\theequation{\thesection.\@arabic\c@equation}
\makeatother
\renewcommand{\theequation}{\arabic{section}.\arabic{equation}}
\newcommand{\R}{\mathbb R}

\newcommand{\N}{\mathbb N}

\def \O{\Omega}
\everymath{\displaystyle}
\newtheorem{theorem}{Theorem}[section]
\newtheorem{pro}{Proposition}[section]

\newtheorem{Def}{Definition}[section]

\newtheorem{cl}{Claim}[section]

\renewcommand{\theequation}{\thesection.\arabic{equation}}
\renewcommand{\thesection}{\arabic{section}}
\renewcommand{\theequation}{\thesection.\arabic{equation}}
\let\ssection=\section\renewcommand{\section}{\setcounter{equation}{0}\ssection}
\begin{document}
\begin{frontmatter}
\title{Multiplicity of solutions for a class of fractional $p(x,\cdot)$-Kirchhoff type \\
problems without the Ambrosetti-Rabinowitz condition}
\author[mk0,mk1,mk2]{M.K. Hamdani}
\ead{hamdanikarim42@gmail.com}
\author[jz1,jz2,jz3]{J. Zuo}
\ead{zuojiabin88@163.com}
\author[ntc1]{N.T. Chung}
\ead{ntchung82@yahoo.com}
\author[dr1,dr2,dr3]{D.D.  Repov\v{s}\corref{cor1}}
\ead{dusan.repovs@guest.arnes.si}
\begin{center}
\address[mk0] {Science and Technology for Defense Laboratory, Military Research Center, Aouina, Tunisia.}
\address[mk1]{Military Aeronautical Specialities School, Sfax, Tunisia.}
\address[mk2]{Mathematics Department, Faculty of Science, University of Sfax, Sfax, Tunisia.}
\address[jz1] {College of Science, Hohai University,  Nanjing  210098, P. R. China.}
\address[jz2] {Faculty of Applied Sciences, Jilin Engineering Normal University, Changchun 130052, P. R. China.}
\address[jz3] {Departamento de Matem\'{a}tica, Universidade Estadual de Campinas, IMECC, Rua S\'{e}rgio Buarque de Holanda, 651, Campinas, SP CEP 13083-859, Brazil.}
\address[ntc1]{Department of Mathematics, Quang Binh University, 312 Ly Thuong Kiet, Dong Hoi, Quang Binh, Vietnam.}
\address[dr1]{Faculty of Education, University of Ljubljana, Kardeljeva pl. 16, SI-1000 Ljubljana, Slovenia.}
\address[dr2]{Faculty of Mathematics and Physics, University of Ljubljana, Jadranska 21, SI-1000 Ljubljana, Slovenia.}
\address[dr3]{Institute of Mathematics, Physics and Mechanics,  Jadranska 19, SI-1000 Ljubljana, Slovenia.}
\cortext[cor1]{Corresponding Author}
\end{center}
\begin{abstract}
We are interested in the existence of solutions for the following fractional $p(x,\cdot)$-Kirchhoff type problem
$$
\left\{\begin{array}{ll}
M \, \left(\displaystyle\int_{\Omega\times \Omega} \ \displaystyle{\frac{|u(x)-u(y)|^{p(x,y)}}{p(x,y) \ |x-y|^{N+p(x,y)s}}} \ dx \, dy\right)(-\Delta)^{s}_{p(x,\cdot)}u = f(x,u), \quad x\in \Omega,  \\    \\
u= 0, \quad x\in \partial\Omega,
\end{array}\right.
$$
where $\Omega\subset\mathbb{R}^{N}$, $N\geq 2$ is a bounded smooth domain, $s\in(0,1),$ 
 $p: \overline{\Omega}\times \overline{\Omega} \rightarrow (1, \infty)$, $(-\Delta)^{s}_{p(x,\cdot)}$ denotes the $p(x,\cdot)$-fractional Laplace operator, $M: [0,\infty) \to [0, \infty),$ and $f: \Omega \times \mathbb{R} \to \mathbb{R}$ are continuous functions.
Using variational methods, especially the symmetric mountain pass theorem due to Bartolo-Benci-Fortunato (Nonlinear Anal. 7:9 (1983), 981-1012), we establish the existence of infinitely many solutions for this problem 
without assuming the Ambrosetti-Rabinowitz condition.
Our main result in several directions extends previous ones which have recently appeared in the literature.
\end{abstract}
\date{}
\end{frontmatter}

\section{Introduction}\label{sec1}

Let $\Omega$ be a smooth bounded domain in $\mathbb{R}^{N}$, $N\geq 2$. Let us consider the following fractional $p(x,\cdot)$-Kirchhoff type problem
\begin{align}\label{eqt1}
\left\{\begin{array}{ll}
M \, \left(\displaystyle\int_{\Omega\times \Omega} \ \displaystyle{\frac{|u(x)-u(y)|^{p(x,y)}}{p(x,y) \ |x-y|^{N+p(x,y)s}}} \ dx \, dy\right)(-\Delta)^{s}_{p(x,\cdot)}u = f(x,u), \quad x\in \Omega, \\   \\
u= 0, \quad x\in \partial\Omega,
\end{array}\right.
 \end{align}
 where $0<s<1$, $p: \overline{\Omega}\times \overline{\Omega} \rightarrow (1, \infty)$ is a continuous function with $sp(x,y)< N$ for all $(x,y) \in \overline{\Omega}\times \overline{\Omega}$, and $M$,$f$ are continuous functions satisfying certain growth conditions to be stated later on. 
 
 The fractional $p(x,\cdot)$-Laplacian operator $(-\Delta)^s_{p(x,\cdot)}$ is, up to normalization factors by the Riesz potential, defined as follows: for each $x\in \Omega$,
\begin{equation}
(-\Delta)^s_{p(x,\cdot)} \varphi(x) = {\rm{p.v.}} 
\displaystyle
\int_{\Omega}\frac{|\varphi(x) - \varphi(y)|^{p(x,y)-2} \ (\varphi(x) - \varphi(y))}{|x-y|^{N+sp(x,y)}} \, dy,
\end{equation}
along any $\varphi\in C^\infty_0(\O)$, where {\sl p.v.} is 
the commonly used abbreviation for the {\sl principal value}.

Throughout this paper, we shall assume that $M: \mathbb R_0^+ := [0, +\infty) \to \mathbb R_0^+$ is a continuous function satisfying the following conditions: \\
$(M_{1}):$ \ \ there exist $\tau_0 > 0$ and  $\gamma\in\left(1,{(p^*_s)^{-}}/{p^{+}}\right)$ such that
\begin{equation*}
tM(t)\leq\gamma\widehat{M}(t),~~~\mbox{for all} ~ t\geq\tau_0,
\end{equation*}
where 
\begin{equation*}
\widehat{M}(t)=\int_{0}^{t}M(\tau) \, d\tau
\end{equation*}
and $p^+$ and $p^-$ 
will be defined in Section 2;\\
$(M_{2}):$ \ \ for every
 $\tau>0$  there exists $\kappa=\kappa(\tau)>0$ such that
\begin{equation*}
M(t)\geq \kappa, ~~~\mbox{for all}~ t\geq\tau.
\end{equation*}

Obviously, the conditions $(M_1)$ and $(M_2)$ are fulfilled for
 the model case: 
 \begin{equation}
 M(t) = a + b \, \theta \, t^{\theta-1}, 
 \ \ \mbox{ where} \ \  
 a, b \geq 0 
 \ \ \mbox{and} \ \ 
 a + b > 0.
 \end{equation} 

 It is worth pointing out that condition $(M_2)$ was originally used to establish multiplicity of solutions for a class of higher order $p(x)$-Kirchhoff equations~\cite{CP}.

In recent years, a lot of attention has been given to problems involving fractional and nonlocal operators. This type of operators arises in
a natural way in many different applications, e.g., image processing, quantum mechanics, elastic mechanics, electrorheological fluids (see \cite{3,4,5,R3} and the references therein).

 In their pioneering paper, Bahrouni and R\u{a}dulescu \cite{BR18}
studied  qualitative properties of the fractional Sobolev space $W^{s,q(x),p(x,y)}(\Omega)$, where $\Omega$ is a smooth bounded domain. Their results have been applied in the variational analysis of a class of nonlocal fractional problems with several variable exponents. 

Recently, by means of approximation and energy methods, Zhang and Zhang \cite{ZZ} have established the existence and uniqueness of nonnegative renormalized solutions for such problems. When $s=1$, the operator degrades to integer order. It has been extensively studied in the literature, see for example 
\cite{
Chung13, 
CN10,  
H2,
HCR,
 HR} and the references therein. In particular, when $p(x,\cdot)$ is a constant, this operator is reduced to the classical fractional $p$-Laplacian operator. 
 
 For studies concerning this operator, we  refer to \cite{
 P8, 
 P9,
  004,
 ZUO1,  005}. We emphasize that, unless the functions $p(x,.)$ and $q(x)$ are constants, the space $W^{s,q(x),p(x,y)}(\Omega)$ does not coincide with the Sobolev space $W^{s,p(x)}(\Omega)$ when $s$ is a natural number, see \cite{5, HR16}. However, because of various applications in physics as well as mathematical finance, the study of nonlocal problems in such spaces is still very interesting.

On the other hand, a lot of interest has in recent years  been devoted
to the study of Kirchhoff-type problems. More precisely,  in 1883 Kirchhoff \cite{kkk}
 established a model given by the following
 equation
\begin{equation}
\rho\frac{\partial ^2u}{\partial t^2}
-\left(\frac{p_0}{\lambda}+\frac{E}{2L}\int_0^L\left|\frac{\partial u}{\partial x}\right|^2dx
\right)\frac{\partial ^2u}{\partial x^2}=0,
\end{equation}
a generalization of the well-known D'Alembert wave equation for free vibrations of elastic strings, where $\rho$, $p_0$, $\lambda$, $E$, $L$ are constants which represent some physical meanings, respectively.

In the study of problem~\eqref{eqt1}, the following Ambrosetti-Rabinowitz condition given in \cite{AR} has been widely used:

(AR): \ \  There exists a constant $\mu>p^{+}$ such that
\begin{equation*}
t  f(x,t) \geq \mu F(x,t) > 0,
\ \ 
\mbox{where}
\ \ 
F(x,t)=\int^{t}_{0}f(x,s) \, ds.
\end{equation*}

Clearly, if the $(AR)$ condition holds, then
\begin{eqnarray}\label{nn}
F(x,t)\geq c_1 \, |t|^\mu-c_2,\;\mbox{ for all } (x,t)\in \Omega\times \R,
\end{eqnarray}
where $c_1$, $c_2$ are two positive constants.

 It is well known that $(AR)$ condition is very important for ensuring the boundedness of the Palais-Smale sequence. When the nonlinear term $f$ satisfies the $(AR)$ condition, many results have been obtained  by using the critical point theory and variational methods, see for example 
 \cite{AHKC19, AB, 2, AOB, DM, DS1, DS2, P3, Ham, MBTZ, BRS, P12, P13}. 
 In particular, 
 Ali {\sl et al.} \cite{AHKC19}
  and 
Azroul {\sl et al.} \cite{2}   have established
 the existence of nontrivial weak solutions for a class of fractional $p(x,\cdot)$-Kirchhoff type problems by using the mountain pass theorem of Ambrosetti and Rabinowitz, direct variational approach, and Ekeland's variational principle.

Since the $(AR)$ condition implies
condition
 \eqref{nn}, one cannot deal with problem~\eqref{eqt1} by 
 using the mountain pass theorem directly if $f(x, t)$ is  $p^+$-asymptotically linear at $\infty$, i.e.
\begin{equation}
\lim\limits_{|t| \to \infty} \dfrac{f(x, t)}{|t|^{p^+-1 }}= l, \quad \mbox{ uniformly in } x\in\Omega,
\end{equation}
where $l$ is a constant. For this reason, in recent years some authors have studied 
problem~\eqref{eqt1} by trying to omit the condition ${(AR)}$, see for example \cite{H2, HHS, LRDZ}. 

Not having the $(AR)$ condition brings great difficulties, so it is natural to
consider if this kind of fractional problems have corresponding results even if the nonlinearity does not satisfy the $(AR)$ condition. 
In fact,  in the absence of Kirchhoff's interference,
Lee at al. \cite{1}
have obtained
 infinitely many solutions to a fractional $p(x)$-Laplacian equation without assuming the $(AR)$ condition, by using the fountain theorem and the dual fountain theorem.

 Inspired by the above works, we consider in this paper 
 the
  fractional $p(x,\cdot)$-Kirchhoff type 
 problem
 without the $(AR)$ condition. 
 Our situation is different from \cite{AHKC19,2} since our Kirchhoff function $M$ belongs to a 
 larger class of functions, whereas the nonlinear term $f$ is $p^+$-asymptotically linear at $\infty$. 
 
 More precisely, let us assume that $f$ satisfies the following global conditions:
\begin{enumerate}
\item[$(F_{1}):$] $f: \Omega \times \mathbb{R}\to \mathbb{R}$ is continuous with $F(x,t) \geq 0,$ for all $(x,t) \in \Omega \times \mathbb{R}$, where 
\begin{equation*}
F(x,t)=\int^{t}_{0}f(x,s) \, ds;
\end{equation*}
\item[$(F_{2}):$]  there exist a function $\alpha \in C(\overline{\Omega})$, $p^{+}<\alpha^{-} \leq \alpha(x)<p_{s}^{*}(x),$ for all $x\in\Omega$, and a number $\Lambda_0>0$ such that for each $\lambda \in (0, \Lambda_0),$  $\epsilon>0$,  there exists $C_\epsilon>0$ such that
\begin{equation*}
f(x,t)\leq (\lambda +\epsilon) \ |t|^{\overline p(x)-1}+ C_\epsilon \,  |t|^{\alpha(x)-1}, ~\mbox{ for all }  (x,t)\in \Omega\times\mathbb{R};
\end{equation*}
\item[$(F_{3}):$] the following is uniformly satisfied on  $ \overline\Omega$
\begin{equation*}
\lim\limits_{|t| \to \infty} \dfrac{F(x, t)}{|t|^{p^+\gamma } }= \infty;
\end{equation*}
\item[$(F_{4}):$] there exist  constants $\mu>p^{+}\gamma$ and $\varpi_{0}>0$ such that
\begin{equation*}
F(x,t) \leq\frac{1}{\mu} \ 
f(x,t) \ 
t 
+
\varpi_{0} \ 
|t|^{p^-}, ~\mbox{ for all }  (x,t)\in \Omega\times\mathbb{R};
\end{equation*}
where $\gamma$ is given by $(M_{1})$;
\item[$(F_{5}):$] the following holds
\begin{equation*}
f(x,-t)=-f(x,t), ~~\mbox{for all}~~ (x,t)\in \Omega\times\mathbb{R}.
\end{equation*}
\end{enumerate}
A simple computation proves that the following function
\begin{eqnarray}\label{yu}
f(x,t)=|t|^{p^+\gamma-2}t \ln^{\alpha(x)}(1+|t|), \ \  \mbox{where} \ \ \alpha(x)> 1,
\end{eqnarray}
does not satisfy ${(AR)}$ condition. However, it is easy to see that $f(x,t)$ in  \eqref{yu} satisfies conditions  $(F_1)-(F_5)$.

 We can now state the definition of (weak) solutions for problem~\eqref{eqt1}   (see Section \ref{sec2}
  for details):
  
\begin{Def}
A function $u\in E_{0}=W_{0}^{s,q(x),p(x,y)}(\Omega)$ is called a (weak) solution of problem~\eqref{eqt1}, if for every $w\in E_{0}$ it satisfies the following
\begin{equation*}
M\Big(\sigma_{p(x,y)}(u)\Big)\int_{\Omega\times \Omega}\displaystyle{\frac{|u(x)-u(y)|^{p(x,y)-2} \ (u(x)-u(y))(w(x)-w(y))}{|x-y|^{N+p(x,y)s}}} \ dx \, dy - \int_\Omega f(x,u) \, w \, dx=0
\end{equation*}
\end{Def}
where
\begin{equation*}
\sigma_{p(x,y)}(u)=\int_{\Omega\times \Omega}\frac{1}{p(x,y)}\frac{|u(x)-u(y)|^{p(x,y)}}{|x-y|^{N+s{p(x,y)}}} \,dx\,dy.
\end{equation*}

The main result of our paper is the following theorem:

\begin{theorem}\label{th1}
Let $q(x),\,p(x,y)$ be continuous variable functions such that $sp(x,y)<N$, $p(x,y) = p(y,x)$ for all $(x,y)\in \overline{\Omega}\times\overline{\Omega}$ and $q(x) \geq p(x,x)$ for all $x \in \overline{\Omega}$. Assume that $f:\Omega \times \mathbb{R}\to \mathbb{R}$ satisfies conditions $(F_1)-(F_5)$ and
that $M: \mathbb R_0^+\to\mathbb R_0^+$ is a continuous function satisfying conditions $(M_1)$ and $(M_2)$. 
Then there exists $\Lambda>0$ such that for each $\lambda \in (0,\Lambda)$, problem~\eqref{eqt1} has a sequence $\{u_n\}_n$ of nontrivial solutions.
\end{theorem}

The paper is organized as follows. In Section \ref{sec2}, we shall introduce the necessary properties of variable exponent Lebesgue spaces and  fractional Sobolev spaces with variable exponent. In Section \ref{sec3}, we shall
 verify the Cerami compactness condition. Finally, in Section \ref{sec4}, we  shall prove Theorem \ref{th1} by means of a version of the  mountain pass theorem. 
\section{Fractional Sobolev spaces with variable exponent}\label{sec2}

For a smooth bounded domain $\Omega$ in $\mathbb{R}^N$, we consider a continuous function $p: \overline{\Omega}\times \overline{\Omega} \rightarrow (1, \infty)$. We assume that $p$ is symmetric, that is,
 \begin{equation*}
 p(x,y)=p(y,x), \ \mbox{ for all} \  (x,y)\in\overline{\Omega}\times \overline{\Omega}
 \end{equation*} 
 and
\begin{equation*}
1< p^-:= \underset{(x,y)\in\overline{\Omega}\times \overline{\Omega}}{\min}p(x,y)\leq p(x,y)\leq p^+ =\underset{(x,y)\in\overline{\Omega}\times \overline{\Omega}}{\max}p(x,y)<\infty.
\end{equation*}
We  also introduce a continuous function $q: \overline{\Omega} \rightarrow \mathbb{R}$ such that
\begin{equation*}
1< q^- := \underset{x\in\overline{\Omega}}{\min} \,q(x) \leq q(x) \leq q^+ :=\underset{x\in\overline{\Omega}}{\max} \,q(x)<\infty.
\end{equation*}
We first give some basic properties of  variable exponent Lebesgue spaces. Set
\begin{equation*}
C_{+}(\overline{\Omega})=\bigg\{r\in C(\overline{\Omega}):1<r(x) ~~\text{for~all}~~ x\in \overline{\Omega}\bigg\}.
\end{equation*}
Given $r\in C_{+}(\overline{\Omega}), $ we define the variable exponent Lebesgue space as
\begin{equation*}
L^{r(x)} (\Omega)=\Bigg\{u: ~ \Omega\rightarrow \mathbb{R}~~ \text{is measurable:} ~\int_{\Omega}|u(x)|^{r(x)} \ dx<\infty\Bigg\},
\end{equation*}
and this  space is endowed with the Luxemburg norm,
\begin{equation*}
|u|_{r(x)}=\inf\Bigg\{\mu>0:\int_{\Omega}\Bigg|\frac{u(x)}{\mu}\Bigg|^{r(x)}dx\leq1\Bigg\}.
\end{equation*}

Then $(L^{r(x)}(\Omega), |\cdot|_{r(x)})$ is a separable reflexive Banach space, see \cite[Theorem $2.5$ and Corollaries $2.7$ and $2.12$]{39}.

Let $\widetilde{r}\in C_{+}(\overline{\Omega})$ be the conjugate exponent of $q$, that is
\begin{equation*}
\frac{1}{r(x)}+\frac{1}{\widetilde{r}(x)}=1\qquad\mbox{for all }x\in\overline{\Omega}.
\end{equation*}
We shall need the following H\"{o}lder inequality, whose proof can be found in \cite[Theorem $2.1$]{39}.
Assume that $v\in L^{r(x)}(\Omega)$ and $u\in L^{\widetilde{r}(x)}(\Omega)$. Then
\begin{equation*}
\Bigg|\int_{\Omega}uv \ dx
\,
\Bigg|
\leq
\Bigg(
\frac{1}{r^{-}}+\frac{1}{\widetilde{r}^{-}}
 \Bigg)
  \ 
 |u|_{r(x)} \ |v|_{\widetilde{r}(\cdot)}\leq 2 \ |u|_{r(x)} \ |v|_{\widetilde{r}(x)}.
\end{equation*}
A modular of the $L^{r(x)}(\Omega)$ space is defined by
\begin{equation*}
\varrho_{r(x)}:L^{r(x)}(\Omega)\rightarrow\mathbb{R}, \ \ u\mapsto\varrho_{r(x)}(u)=\int_{\Omega}|u(x)|^{r(x)} \ dx.
\end{equation*}
Assume that $u\in L^{r(x)}(\Omega)$ and $\{u_n\}\subset L^{r(x)}(\Omega)$.
Then the following assertions hold (see \cite{5}):
\begin{align*}
\noindent(1)~~&|u|_{r(x)}<1 \ \  (\mbox{resp.,}=1,>1)\Leftrightarrow \varrho_{r(x)}(u)<1 \ \  (\mbox{resp.,}=1,>1),\\
\noindent(2)~~& |u|_{r(x)}<1\Rightarrow |u|_{r(x)}^{r^{+}}\leq\varrho_{r(x)}(u)\leq |u|_{r(x)}^{r^{-}},\\
\noindent(3)~~& |u|_{r(x)}>1\Rightarrow |u|_{r(x)}^{r^{-}}\leq\varrho_{r(x)}(u)\leq |u|_{r(x)}^{r^{+}},\\
\noindent(4)~~&\lim\limits_{n\rightarrow\infty}|u_{n}|_{r(x)}=0 
\ \ (\mbox{resp.,} \ \
= \infty)\Leftrightarrow\lim\limits_{n\rightarrow\infty}\varrho_{r(x)}(u_{n})=0
 \ \ 
 (\mbox{resp.,} \ \ 
= \infty),\\
\noindent(5)~~&\lim\limits_{n\rightarrow\infty}|u_{n}-u|_{r(x)}=0\Leftrightarrow\lim\limits_{n\rightarrow\infty}\varrho_{r(x)}(u_{n}-u)=0.
\end{align*}

Given $s\in (0,1)$ and the functions $p(x,y)$, $q(x)$ as we mentioned above, the fractional Sobolev space with variable exponents via the Gagliardo approach $E=W^{s,q(x),p(x,y)}(\Omega)$ is  defined as follows
\begin{equation*}
E=\Bigg{\{}u\in L^{q(x)}(\Omega): \int_{\Omega\times\Omega}  \  
\frac{|u(x)-u(y)|^{p(x,y)}}{\mu^{p(x,y)}
 \ 
|x-y|^{N+sp(x,y)}} \ dx \, dy<\infty,\,\hbox{ \text{for some }} \mu >0\Bigg{\}}.
\end{equation*}

Let
\begin{equation*}
[u]_{s,p(x,y)}= \inf \Bigg{\{}\mu >0:\,\int_{\Omega\times\Omega} \ 
 \frac{|u(x)-u(y)|^{p(x,y)}}{\mu^{p(x,y)}
  \ 
 |x-y|^{N+sp(x,y)}} \ dx \, dy<1\Bigg{\}},
\end{equation*}
be the variable exponent Gagliardo seminorm and define
\begin{equation*}
\|u\|_E= [u]_{s,p(x,y)}+|u|_{q(x)}.
\end{equation*}
Then $E$ equipped with the norm $\|\cdot\|_E$ becomes
a Banach space.

\begin{pro}\label{lem2}The following properties hold:
\begin{enumerate}
\item[(1)] If $1\leq [u]_{s,p(x,y)}<\infty$, then
\begin{equation*}
([u]_{s,p(x,y)})^{p_-}\leq \int_{\Omega\times\Omega} \ 
 \frac{|u(x)-u(y)|^{p(x,y)}} 
 {|x-y|^{N+sp(x,y)}} \ dx \, dy\leq ([u]_{s,p(x,y)})^{p^+}.
\end{equation*}
\item[(2)] If $[u]_{s,p(x,y)}\leq 1$, then
\begin{equation*}
([u]_{s,p(x,y)})^{p^+}\leq \int_{\Omega\times\Omega} \ 
\frac{|u(x)-u(y)|^{p(x,y)}}  
{|x-y|^{N+sp(x,y)}} \ dx \, dy\leq ([u]_{s,p(x,y)})^{p_-}.
\end{equation*}
\end{enumerate}
\end{pro}

\noindent
Given
 $u\in W^{s,q(x),p(x,y)}(\Omega)$, we set
\begin{equation*}\label{e1}
\rho(u)=\int_{\Omega\times\Omega} \ 
\frac{|u(x)-u(y)|^{p(x,y)}}  
{|x-y|^{N+sp(x,y)}} \ dx \, dy + \int_\Omega |u|^{q(x)} dx
\end{equation*}
and
\begin{equation*}\label{e2}
\|u\|_{\rho}= \inf \Bigg{\{} \mu >0:\;\rho\bigg(\frac{u}{\mu}\bigg) \leq 1\Bigg{\}}.
\end{equation*}

It is well-known that $\|\cdot\|_{\rho}$ is a norm which is equivalent to the norm $\|\cdot\|_{W^{s,q(x),p(x,y)}(\Omega)}.$
By Lemma 2.2 in \cite{ZZ}, $(W^{s,q(x),p(x,y)}(\Omega), \,\|\cdot\|_{\rho})$ is uniformly convex and $W^{s,q(x),p(x,y)}(\Omega)$ is a reflexive Banach space.

We denote our workspace 
$E_0= W_0^{s,q(x),p(x,y)}(\Omega),$
the closure of 
$C_0^{\infty}(\Omega)$
 in $E$. Then $E_0$ is a 
reflexive
Banach space with the norm 
\begin{equation*}
\|\cdot\|_{E_0}=[u]_{s,p(x,y)}.
\end{equation*}

A thorough variational analysis of the problems with variable exponents has been developed in
the  monograph by R\u{a}dulescu and Repov\v{s} \cite{RD}. The following result provides a compact embedding  into variable exponent Lebesgue spaces.

\begin{theorem}[see \cite{ZZ}]\label{l1}
Let  $\Omega \subset \mathbb{R}^n$ be a smooth bounded domain and $s\in (0,1)$. Let $q(x),\,p(x,y)$ be continuous variable exponents
such that 
\begin{equation*}
sp(x,y)<N, \ \ \mbox{for} \ \  (x,y)\in \overline{\Omega}\times\overline{\Omega}
\ \ \mbox{and} \ \ 
 q(x) \geq p(x,x),
 \ \ \mbox{for all} \ \ 
 x \in \overline{\Omega}.
 \end{equation*} 
Assume that $\tau: \overline{\Omega}\longrightarrow(1,\infty)$ is a continuous function such that
\begin{equation*}
p^*(x)= \frac{Np(x,x)}{N-sp(x,x)} > \tau(x) \geq \tau^->1,\;\;\mbox{for all}\;\; x\in \overline{\Omega}.
\end{equation*}
Then there exists a constant $C=C(N,s,p,q,r,\Omega)$ such that for every $u\in W^{s,q(x),p(x,y)}$,
\begin{equation}\label{sobolev}
|u|_{\tau(x)} \leq C\|u\|_{E}.
\end{equation}

That is, the space $W^{s,q(x),p(x,y)}(\Omega)$ is continuously embeddable in $L^{\tau(x)}(\Omega)$. Moreover, this embedding is compact. In addition, if $u\in W_0^{s,q(x),p(x,y)}$, the following inequality holds
\begin{equation*}
|u|_{\tau(x)}\leq C\|u\|_{E_{0}}.
\end{equation*}
\end{theorem}

\begin{theorem}[see \cite{BR18}]\label{111}
For all $u, v\in E_0$, we consider the following 
operator
 $I:E_0\to {E^{*}_0}$ such that
\begin{equation*}
 \langle I(u),v\rangle =\int_{\Omega\times\Omega} \ 
  \frac{|u(x)-u(y)|^{p(x,y)-2} \ 
  (u(x)-u(y))(v(x)-v(y))} {|x-y|^{N+sp(x,y)}} \ dx \, dy. 
 \end{equation*}
Then the following properties hold:
\begin{enumerate}
\item[(1)] $I$ is a bounded and strictly monotone operator.
\item[(2)] $I$ is a mapping of type $(S_+)$, that is, 
\begin{equation*}
\mbox{if} \ \ 
u_n\rightharpoonup u \in E_0
\ \ \mbox{and} \ \ 
 \displaystyle\limsup_{n\rightarrow \infty} I(u_n)(u_n-u)\leq 0,
 \ \ \mbox{then} \ \ 
 u_n\rightarrow u \in E_0.
 \end{equation*}
 \item[(3)] $I:E_0\to {E^{*}_0}$ is a homeomorphism.
\end{enumerate}
\end{theorem} 
\section{The Cerami Compactness Condition}\label{sec3}

Let us consider the Euler-Lagrange functional associated to problem~\eqref{eqt1}, defined by $J_\lambda: E_0 \to \mathbb{R}$
\begin{eqnarray}\label{func}
J_\lambda(u)&=&\widehat{M}\Big(\sigma_{p(x,y)}(u)\Big)-\int_{\Omega}F(x,u) \ dx.
 \end{eqnarray}
Note that $J_\lambda$ is a $C^1(E_0, \R)$ functional and
\begin{eqnarray} \label{deriv-fun}
\langle J_\lambda'(u), {w}\rangle &=M\Big(\sigma_{p(x,y)}(u)\Big)\int_{\Omega\times \Omega}
 \ 
 \displaystyle{\frac{|u(x)-u(y)|^{p(x,y)-2} \ 
 (u(x)-u(y))(w(x)-w(y))}{|x-y|^{N+p(x,y)s}}} \ dx \, dy \nonumber\\
&-\int_\Omega f(x,u) \, w \ dx
\end{eqnarray}
for all $w\in E_0$. Therefore
 critical points of $J_\lambda$ are weak solutions of 
 problem \eqref{eqt1}.\\

In order to prove our main result (Theorem \ref{th1}), we  recall the definition of the
Cerami compactness condition \cite{PRR}.

\begin{Def}\label{def1}
We say that $J_\lambda$ satisfies the Cerami compactness condition at the level $c \in \R$ ((Ce)$_c$ condition for short), if 
every sequence $\{u_n\}_n\subset E_0$, i.e., $J_\lambda(u_n) \to c$ and 
\begin{equation*}
\|J_\lambda'(u_n)\|_{{E_0}^*}(1+\|u_n\|_{E_0})\to 0,
\ \ \mbox{ as} \ \  n\to \infty,
\end{equation*} 
admits a strongly convergent subsequence in $E_0$. If $J_\lambda$ satisfies the (Ce)$_c$ condition for any $c \in \R$ then we say that $J_\lambda$ satisfies the Cerami compactness condition.
\end{Def}

\begin{cl}\label{Bounded}
 Under assumptions of Theorem \ref{th1}, every $(Ce)_c$ sequence $ \{u_n\}_n\subset E_0$ of $J_\lambda$ is bounded in $E_0$.
\end{cl}

{\bf Proof.}
 Let $\{u_n\}_n$ be a $(Ce)_c$-sequence of $J_\lambda$. Then
 \begin{equation}\label{cond7}
J_\lambda(u_n)\to  c \ \  \ \mbox{and}  \ \ \ 
\|J_\lambda'(u_n)\|_{{E_0}^*}(1+\|u_n\|_{E_0})\to 0.
\end{equation}
First, we prove that the sequence $\{u_n\}_n$ is bounded in $E_0$. To this end, we argue by contradiction. So suppose that $\|u_n\|_{E_0} \to\infty$, as $n\to \infty$. We define the sequence 
$\{v_n\}_n$ by
\begin{equation*}
v_n=\frac{u_n}{\|u_n\|_{E_0}}, \ \  n \in \N.
\end{equation*}
 It is clear that $\{v_n\}_n\subset E_0$  and $\|v_n\|_{E_0} = 1$ for all $n\in \N$.  Passing, if necessary, to a subsequence, we may assume that
 \begin{eqnarray}\label{cond cvg}
   v_n&\rightharpoonup& v\quad \text{in } \ \ E_0,\nonumber\\
   v_n&\to&  v \quad \text{in } \ \ L^{\tau(x)}(\O),\; 1\leq \tau(x)<p^*(x),\\
   v_n(x)&\to&  v(x) \quad \text{a.e. on } \O.\nonumber
 \end{eqnarray}
Let $\O_\natural:=\{x \in \O : v(x)\neq 0\}$. If $x \in\O_\natural$, then it follows from \eqref{cond cvg} that
\begin{equation*}
\lim\limits_{n\to
\infty}v_n(x)=\lim\limits_{n\to
\infty}\frac{u_n}{\|u_n\|_{E_0}}=v(x)\neq 0.
\end{equation*}
This means that
\begin{equation*}
|u_n(x)|=|v_n(x)| \ 
\|u_n\|_{E_0}\to +\infty \  \mbox{ a.e. on} \ \O_\natural, \ \  \mbox{as}  \ n\to
\infty.
\end{equation*}
Moreover, it follows by condition $(F_3)$ and Fatou's Lemma that for each $x\in \O_\natural$,
\begin{eqnarray}\label{cond10}
+\infty=\lim\limits_{n\to
\infty}\int_{\Omega}\frac{|F(x,u_n(x))|}{|u_n(x)|^{p_+\gamma}}
\,
\frac{|u_n(x)|^{p_+\gamma}}{\|u_n(x)\|^{p_+\gamma}_{E_0}}
 \ 
  dx=\lim\limits_{n\to
\infty}\int_{\Omega}\frac{|F(x,u_n(x))|
\,
|v_n(x)|^{p_+\gamma}}
{|u_n(x)|^{p_+\gamma}} \ dx.
\end{eqnarray}
Condition $(M_1)$ gives
\begin{equation}\label{aM}
\widehat{M}(t) \leq \widehat{M}(1)
 \ 
 t^\gamma, \quad\mbox{ for all } t\geq 1.
\end{equation}
Now, since $\|u_n\|_{E_0}>1$, it follows by \eqref{func}, \eqref{cond7} and \eqref{aM} that
\begin{eqnarray*}\label{cond11}
\int_{\Omega}F(x,u_n) 
\ 
dx&\leq&\widehat{M}\Big(\sigma_{p(x,y)u_n}\Big)+C\nonumber\\
&\leq&\frac{\widehat{M}(1)}{(p^-)^\gamma}
 \
 (\sigma_{p(x,y)}(u_n))^\gamma+C\nonumber\\
&\leq& \frac{\widehat{M}(1)}{(p^-)^\gamma} 
\ 
\|u_n\|_{E_0}^{\gamma p^+}+C,
\end{eqnarray*}
for all $n\in \N$. We can now conclude that
\begin{eqnarray}\label{contr}
\lim\limits_{n\to \infty}\int_{\Omega}\frac{F(x,u_n)}{\|u_n\|_{E_0}^{p_+\gamma}}
 \
 dx&\leq& \lim\limits_{n\to \infty}\left(\frac{\widehat{M}(1)}{(p^-)^\gamma}+\frac{C}{\|u_n\|_{E_0}^{\gamma p^+}}\right).
\end{eqnarray}
From \eqref{cond10} and \eqref{contr} we obtain
\begin{equation*}
+\infty\leq \frac{\widehat{M}(1)}{(p^-)^\gamma},
\end{equation*}
which is a contradiction. Therefore 
\begin{equation*}
|\O_\natural| = 0
\ \mbox{ and} \ \
  v(x) = 0
  \ \ 
  \mbox{ a.e. on} \ \ 
    \O.
\end{equation*}
It follows from $(M_{1})$, $(M_{2})$,  $(F_{4}),$ and since $v_{n}\rightarrow v=0$ in $L^{p^-}(\O)$, that
\begin{align*}
&  \frac{1}{\|u_{n}\|_{E_{0}}^{p^{-}}}
\
\left(J_\lambda(u_{n})-\frac{1}{\mu}J_\lambda^{'}
 (u_{n})
 \
 u_{n}\right) \\
\geq &
\
\frac{1}{\|u_{n}\|_{E_{0}}^{p^{-}}}
\
\Bigg[\widehat{M}\Big(\sigma_{p(x,y)}(u_n)\Big)-\int_{\Omega}F(x,u_n) 
\
dx\\
&-\frac{1}{\mu}
\
M\Big(\sigma_{p(x,y)}(u_n)\Big)\int_{\Omega\times \Omega}\displaystyle{\frac{|u_n(x)-u_n(y)|^{p(x,y)}}{|x-y|^{N+p(x,y)s}}} \ dx \, dy +\frac{1}{\mu}\int_\Omega f(x,u_n)
\
u_n
\
dx\Bigg]\\
\geq&
\
\frac{1}{\|u_{n}\|_{E_{0}}^{p^{-}}}
\
\Bigg[\frac{1}{\gamma}
\
M
\Big(
\sigma_{p(x,y)}(u_n)\Big)
\
\sigma_{p(x,y)}(u_n)-\frac{1}{\mu}
\
M\Big(
\sigma_{p(x,y)}(u_n)\Big)\int_{\Omega\times \Omega}
\
\displaystyle{\frac{|u_n(x)-u_n(y)|^{p(x,y)}}{|x-y|^{N+p(x,y)s}}}
 \ dx \, dy\\
&-\varpi_{0}\int_{\Omega}|u_n|^{p^-}dx\Bigg]\\
\geq&
\
\frac{1}{\|u_{n}\|_{E_{0}}^{p^{-}}}
\
\Bigg[\bigg(\frac{1}{\gamma p^{+}}-\frac{1}{\mu}\bigg)
\
M\Big(\sigma_{p(x,y)}(u_n)\Big)\int_{\Omega\times \Omega}\displaystyle{\frac{|u_n(x)-u_n(y)|^{p(x,y)}}{|x-y|^{N+p(x,y)s}}} \ dx \, dy
-
\varpi_{0}\int_{\Omega}|u_n|^{p^-} \ 
dx\Bigg]\\
\geq&
\
\bigg(\frac{1}{\gamma p^{+}}-\frac{1}{\mu}\bigg)
\
\kappa-\lambda
\
\varpi_{0}\int_{\Omega}|v_n|^{p^{-}}
 \ 
 dx,
\end{align*}
which means that 
\begin{equation*}
0\geq\bigg(\frac{1}{\gamma p^{+}}-\frac{1}{\mu}\bigg)
\,
\kappa, \ \ \mbox{ as} \ \  n\rightarrow\infty.
\end{equation*}
This is a contradiction.
 As a consequence, we can conclude that Cerami sequence $ \{u_n\}_n$ is
 indeed
 bounded. This completes the proof of
 Claim~\ref{Bounded}.
 \qed \\

We now complete the verification of the Cerami compactness condition $(Ce)_c$ for $J_\lambda$.

\begin{cl}\label{ce cvg}
The functional $J_\lambda$ satisfies condition $(Ce)_c$  in $E_0$.
\end{cl}

{\bf Proof. } Let $\{u_n\}_n$ be a $(Ce)_c$ sequence for $J_\lambda$ in $E_0$. Claim \ref{Bounded} asserts that $\{u_n\}_n$ is bounded in $E_0$. By Theorem \ref{l1}, the embedding $E_0 \hookrightarrow L^{\tau(x)}(\O)$ is compact, where $1\leq \tau(x)<p^*(x)$. Since $E_0$ is a reflexive Banach space, passing, if necessary, to a subsequence, still denoted by $\{u_n\}_n$, there exists $u\in E_0$ such that
\begin{equation}\label{cvg}
u_n\rightharpoonup u\mbox{ in } E_0,\; u_n \to u \mbox{ in } L^{\tau(x)}(\O),\;\; u_n(x)\to u(x), \mbox{ a.e. on } \O.
\end{equation}
By virtue of
 \eqref{cond7},
  we get
\begin{eqnarray}\label{s+}
&\langle J_\lambda'(u_n), {u_n-u}\rangle =\nonumber\\&M\Big(\sigma_{p(x,y)}(u_n)\Big)\int_{\Omega\times \Omega}\displaystyle{\frac{|u_n(x)-u_n(y)|^{p(x,y)-2}
\
(u_n(x)-u_n(y))
\
((u_n(x)-u(x))-(u_n(y)-u(y)))}{|x-y|^{N+p(x,y)s}}}
\ dx \, dy \nonumber\\
& -\int_\Omega f(x,u_n)
\
(u_n-u)
\
dx \to 0.
\end{eqnarray}
 Now, by condition $(F_2)$,
\begin{equation}\label{ss}
 |f(x,u_n)| \leq (\lambda + \epsilon)
 \,
  |u_n| ^{\overline{p}(x)-1}+C_\epsilon
 \,
 |u_n| ^{\alpha(x)-1}.
\end{equation}
It follows from \eqref{cvg}, \eqref{ss} and Proposition~\ref{sobolev}
that
\begin{eqnarray*}
\Bigg|\int_\O f(x,u_n)
\
(u_n-u)
\
dx
\
\Bigg|  
&\leq& \int_\O (\lambda + \epsilon) 
\
|u_n| ^{\overline{p}(x)-1}|u_n-u| \ dx +\int_\O C_\epsilon
\
 |u_n| ^{\alpha(x)-1}
 \
 |u_n-u| 
\
dx\\
&\leq& (\lambda + \epsilon)
\
{\Big|{|u_n|} ^{\overline{p}(x)-1}\Big|}_{\frac{\overline{p}(x)}{\overline{p}(x)-1}}{| u_n-u|}_{\overline{p}(x)} 
   +C_\epsilon {\Big|{|u_n|} ^{\alpha(x)-1}\Big|}_{\frac{\alpha(x)}{\alpha(x)-1}}|u_n-u| _{\alpha(x)}\\
&\leq& (\lambda + \epsilon)\max\left\{\|u_n\|^{p_+-1}_{E_0}, \|u_n\|^{p_--1}_{E_0}\right\}{|u_n-u|} _{p(x)}\\&\qquad+&C_\epsilon\max\left\{\|u_n\|^{\alpha_+-1}_{E_0}, \|u_n\|^{\alpha_--1}_{E_0}\right\}_{\frac{\alpha(x)}{\alpha(x)-1}}|u_n-u| _{\alpha(x)}\\&\to& 0, \ \mbox{as} \  n\to \infty,
\end{eqnarray*}
which implies that
\begin{equation}
\label{cvgf}
\lim_{n\to \infty}\int_{\O}f(x,u_n)
\
(u_n-u)
\
dx=0.
\end{equation}
 Therefore we can infer from \eqref{s+} and \eqref{cvgf} that
\begin{equation*}
M\Big(\sigma_{p(x,y)}(u_n)\Big)\int_{\Omega\times \Omega}\displaystyle{\frac{|u_n(x)-u_n(y)|^{p(x,y)-2}
\
(u_n(x)-u_n(y))
\
((u_n(x)-u(x))-(u_n(y)-u(y)))}{|x-y|^{N+p(x,y)s}}}
\
dx
\,
dy \to 0.
\end{equation*}
Since $\{u_n\}_n$ is bounded in $E_0$, using $(M_2)$,
we can conclude that
 the sequence of positive real numbers $\left\{M\Big(\sigma_{p(x,y)}(u_n)\Big)\right\}$ is bounded from below by some positive number for $n$ large enough. Invoking  Theorem \ref{111}, we can
 deduce that $u_n \to u$ strongly in $E_0$.
 This completes the proof of Claim~\ref{ce cvg}. \qed
\section{Proof of Theorem \ref{th1}}\label{sec4}
To prove Theorem \ref{th1}, we shall use the following 
symmetric
mountain pass theorem.

\begin{theorem}[see \cite{T11T,T11CC}] \label{TH-symm}
 Let  $X=Y\oplus Z$  be an infinite-dimensional Banach space,  where $Y$  is finite-dimensional, and let $I\in C^1(X,\mathbb{R})$.
 Suppose that:
\begin{itemize}
\item[(1)] $I$  satisfies $(Ce)_c$-condition, for all $c>0$;
\item[(2)] $I(0)=0$, $I(-u)=I(u)$,  for all  $u\in X$;
\item[(3)] there exist constants $\rho$,~$a>0$  such that
  $I|_{\partial B\rho\cap Z}$ $\geq a$;
\item[(4)]  for every finite-dimensional subspace
 $\widetilde X\subset X$, there is $R=R(\widetilde X)>0$  such that
 $I(u)\leq0$  on  $\widetilde X\backslash B_R$.
\end{itemize}
Then $I$
 possesses an unbounded sequence of critical values.
\end{theorem}

  Let us first verify that  functional $J_\lambda$ satisfies the mountain pass geometry.

\begin{cl}\label{lemma1}
Under the hypotheses of Theorem \ref{th1}, there exists $\Lambda>0$ such that for each $\lambda \in (0, \Lambda)$, we can choose $\rho > 0$ and $a > 0$ such that
\begin{equation*}
J_\lambda(u) \geq a > 0,
\ \ \mbox{ for all} \ \ 
u\in {E_0}
\ \ \mbox{ with} \ \ 
 \| u\| =\rho.
 \end{equation*}
\end{cl}

{\bf Proof.}  Let $\rho\in(0,1)$ and $u\in E_0$ be such that $\| u\|_{E_0} =\rho$. By assumption $(F_2)$, for every $\epsilon>0$, there exists $C_\epsilon>0$ such that
\begin{equation} \label{condg0}
| F(x,t)| \leq\frac{(\lambda + \epsilon)}{\overline{p}(x)}|t| ^{\overline{p}(x)}+\frac{C_\epsilon}{\alpha(x)}|t| ^{\alpha(x)},\quad\mbox{ for all } x\in\Omega, t\in \mathbb{R}.
\end{equation}
Moreover, $(M_2)$ gives
\begin{equation}\label{condm}
\widehat{M}(t) \geq \widehat{M}(1)
\
t^\gamma, \quad\mbox{ for all } \  t\in [0, 1],
\end{equation}
whereas
 $(M_1)$ implies that $\widehat{M}(1) > 0$. Thus, using \eqref{condg0}, \eqref{condm} and \eqref{sobolev},
we obtain for all $u\in {E_0}$, with $\| u\|_{E_0} =\rho$,
 \begin{eqnarray}\label{inf}
J_\lambda(u)&=&\widehat{M}\Big(\sigma_{p(x,y)}(u)\Big)- \int_{\Omega}F(x,u) 
\
dx \nonumber\\
&\geq&\frac{\widehat{M}(1)}{(p^+)^\gamma}
\
\Big(\sigma_{p(x,y)}(u)\Big)^\gamma - \int_\Omega\frac{\lambda + \epsilon}{\overline{p}(x)}
\
|u| ^{\overline{p}(x)}
\
dx- \int_\Omega\frac{C_\epsilon}{\alpha(x)}
\
|u| ^{\alpha(x)}
\
dx \nonumber\\
&\geq&\frac{\widehat{M}(1)}{(p^+)^\gamma}\min\{ \|u\|_{E_0}^{\gamma p^+},\|u\|_{E_0}^{\gamma p^-}\}-\frac{\epsilon+ \lambda}{p^-}\max\{|u|_{p(x)}^{ p^+},|u|_{p(x)}^{ p^-}\}-\frac{C_\epsilon}{\alpha^-}\max\{|u|_{\alpha(x)}^{\alpha^+},|u|_{\alpha(x)}^{\alpha^-}\}\nonumber\\
&\geq&\frac{\widehat{M}(1)}{(p^+)^\gamma}
\
\|u\|_{E_0}^{\gamma p^+}-c_1(\epsilon+ \lambda)\|u\|_E^{ p^-}-c_2\|u\|_{E_0}^{\alpha^-}\nonumber\\
& \geq & \rho^{\gamma p^+}\left(\frac{\widehat{M}(1)}{(p^+)^\gamma}-c_1(\epsilon+ \lambda)
\
\rho^{ p^--\gamma p^+}-c_2
\
\rho^{\alpha^--\gamma p^+}\right),
 \end{eqnarray}
where $\rho = \|u\|_{E_0}$. Since $\epsilon>0$ is arbitrary, let us choose 
\begin{equation}\label{fd}
\epsilon = \frac{\widehat{M}(1)}{2c_1(p^+)^\gamma}
\
 \rho^{\gamma p^+-p^-} > 0.
\end{equation}
 Then by \eqref{inf}  and \eqref{fd}, we obtain
 \begin{eqnarray}\label{inf1}
J_\lambda(u) \geq \rho^{\gamma p^+}\left(\frac{\widehat{M}(1)}{2(p^+)^\gamma}-\lambda c_1
\
\rho^{ p^--\gamma p^+}-c_2
\
\rho^{\alpha^--\gamma p^+}\right).
 \end{eqnarray}
Now, for each $\lambda > 0,$ we define a continuous function, $g_\lambda : (0, \infty) \to\R$,
\begin{equation*}
g_\lambda(s)=\lambda c_1s^{p^--\gamma p^+}+c_2s^{\alpha^--\gamma p^+}.
\end{equation*}
Since $1<p^-<\gamma p^+<\alpha^-$,  it follows that 
\begin{equation*}
\lim_{s\to 0^+}g_\lambda(s) = \lim_{s\to +\infty}g_\lambda(s) =+\infty.
\end{equation*}
 Thus we can find the  infimum of $g_\lambda$. Note that equating
\begin{equation*}
g'_\lambda(s)=\lambda c_1(p^--\gamma p^+)s^{p^--\gamma p^+-1}+c_2(\alpha^--{\gamma p^+})s^{\alpha^--{\gamma p^+}-1}=0,
\end{equation*}
we get 
\begin{equation*}
s_0=s=\widetilde{C}\lambda^\frac{1}{\alpha^--p^-},
\end{equation*}
 where 
\begin{equation*}
\widetilde{C} := \left(\frac{ c_1(\gamma p^+-p^-)}{c_1(\alpha^--\gamma p^+)}\right)^\frac{1}{\alpha^--p^-} >0.
\end{equation*}
Clearly, $s_0 > 0$. 
It can
also 
 be checked that $g''_\lambda(s_0) > 0$ and hence the
infimum of $g_\lambda(s)$ is achieved at $s_0$. 

Now, observing that
\begin{equation*}
g_\lambda(s_0)= \left(c_1\widetilde{C}^{p^--\gamma p^+}+c_2\widetilde{C}^{\alpha^--\gamma p^+}\right) \lambda^\frac{\alpha^--\gamma p^+}{\alpha^--p^-} \ \rightarrow \
 0,  \ \mbox{ as } \  \lambda\to 0^+,
 \end{equation*}
we can infer from  \eqref{inf1} that there exists $0 < \Lambda < \Lambda_0$ (see $(F_2)$) such that all for all $\lambda\in  \left(0,\Lambda\right)$
we can choose $\rho$ small enough and $\alpha > 0$ such that 
\begin{equation*}
J_\lambda(u) \geq a > 0,
\ \ \mbox{ for all} \ \  
u \in E_0
\ \ \mbox{ with} \ \ 
\|u\|_{E_0}=\rho.
\end{equation*}
This completes the proof of Claim~\ref{lemma1}.
\qed

\begin{cl}\label{lemma2}
Under the hypotheses of Theorem \ref{th1}, for every finite-dimensional subspace $W\subset E_0$ there exists $R= R(W)>0$ such that
\begin{equation*}
J_\lambda(u)\leq 0,
\ \ \mbox{for all} \ \ 
 u\in W,
 \ \ \mbox{with} \ \  \|u\|_{E_0}\geq R.
 \end{equation*}
\end{cl}
{\bf Proof.} In view of $(F_3),$ we know that for all $A>0$, there exists $C_A>0$ such that
\begin{equation}\label{eF}
F(x,t)\geq A
\
|t| ^{\gamma p^+}- C_A,\; \mbox{for all} \ \; (x,u)\in
\O\times \mathbb{R}.
\end{equation}
Again, $(M_2)$ gives
\begin{equation}\label{eM}
\widehat{M}(t) \leq \widehat{M}(1)
\
t^\gamma, \quad\mbox{ for all } \ \ t\geq 1,
\end{equation}
with $\widehat{M}(1) > 0$ by $(M_1)$. By \eqref{eF} and \eqref{eM} we have
\begin{align*}
J_\lambda(u)&=\widehat{M}\Big(\sigma_{p(x,y)}(u)\Big)-\int_{\Omega}F(x,u) 
\
dx\\
&\leq \frac{\widehat{M}(1)}{(p^-)^\gamma}
\
(\sigma_{p(x,y)}(u))^\gamma -  A\int_{\O} |u| ^{\gamma p^+} 
\
dx + C_A | \O|\\
&\leq \frac{\widehat{M}(1)}{(p^-)^\gamma}
\
\|u\|_{E_0}^{\gamma p^+} - A\int_{\O} |u| ^{\gamma p^+}
\
 dx + C_A | \O|.
\end{align*}
Consequently, since $\| u\|_{E_0} > 1$, all norms on the finite-dimensional space $W$  are equivalent, so there is $C_W > 0$ such that
\begin{equation*}
\int_{\O}| u| ^{\gamma p^+}dx\geq C_W
\
\|u\|_{E_0} ^{\gamma p^+}.
\end{equation*}
Let $R= R(W)>0$.
Then  for all $u\in W$ with $\|u\|_{E_0}\geq R$ we obtain
 \begin{equation}\label{last}
  J_\lambda(u) \leq \|u\|_{E_0} ^{p^+\gamma}\left(\frac{\widehat{M}(1)}{(p^-)^\gamma}- AC_W\right)+ C_A | \O|.
  \end{equation}
  
So choosing in  inequality~\ref{last},
\begin{equation*}
A=\frac{2\widehat{M}(1)}{C_W(p^-)^\gamma},
\end{equation*} 
we can
conclude that   
\begin{equation*}
J_\lambda(u)\leq 0,
\ \ \mbox{for all} \ \ 
 u\in W
 \ \ \mbox{with} \ \  \|u\|_{E_0}\geq R.
 \end{equation*}
 
This completes the proof of Claim~\ref{lemma2}.\qed
\subsection*{ \bf  Proof of Theorem \ref{th1}.}
Obviously,
$J_\lambda(0)=0$ and by condition $(F_5)$, $J_\lambda$ is an even functional. Invoking Claims
 \ref{Bounded},
 \ref{ce cvg},
  \ref{lemma1}, and
   \ref{lemma2}, and Theorem \ref{TH-symm},
   we can now conclude
   that
     there indeed exists an unbounded
sequence of solutions of problem \eqref{eqt1}.
This completes the proof of Theorem~\ref{th1}.\qed
\section*{Availability of data and materials}
Not applicable.
\section*{Competing interests}
The authors declare that they have no competing interests.
\section*{Funding}
The first author was supported by the Tunisian Military Research Center for Science and Technology Laboratory. The second author was supported by the Fundamental Research Funds for Central Universities (2019B44914) and the National Key Research and Development Program of China (2018YFC1508100), the China Scholarship Council (201906710004), and the Jilin Province Undergraduate Training Program for Innovation and Entrepreneurship (201910204085). The third author was supported by the Vietnam National Foundation for Science and Technology Development (NAFOSTED) (Grant N.101.02.2017.04). The fourth author was supported by the Slovenian Research Agency grants P1-0292, N1-0114, N1-0083, N1-0064, and J1-8131.
\section*{Authors' contributions}
The authors declare that their contributions are equal.
\section*{Acknowledgements}
The authors would like to thank the referees for their comments and suggestions.


\begin{thebibliography}{777}
\bibliographystyle{alpha}
\bibitem{AHKC19}
K. B. Ali, M. Hsini, K. Kefi and  N. T. Chung,
{\sl On a nonlocal fractional $p(.,.)$-Laplacian problem with competing nonlinearities}, 
Complex Anal. Operator Theory $\mathbf{13}$:3 (2019), 1377-1399.

\bibitem{AB}
C. Alves and G. Molica Bisci, 
{\sl A compact embedding result for anisotropic Sobolev spaces associated to a strip-like domain and some applications}, 
J. Math. Anal. Appl., published online, https://doi.org/10.1016/j.jmaa.2019.123490.

\bibitem{AR}
A. Ambrosetti and P. Rabinowitz,
{\sl  Dual variational methods in critical point theory and applications},
J. Funct. Anal. $\mathbf{14}$ (1973), 349-381.

\bibitem{AOB}
V. Ambrosio, L. D'Onofrio and G. Molica Bisci, 
{\sl Perturbation methods for nonlocal Kirchhoff-type problems},
Fract. Calc. Appl. Anal. {\bf 20} (2017), 829-853.

\bibitem{2}
E. Azroul, A. Benkirane, M. Shimi and M. Srati,
{\sl On a class of fractional $p(x)$-Kirchhoff type problems},
 Appl. Anal., published online, https://www.tandfonline.com/doi/full/10.1080/00036811.2019.1603372.

\bibitem{BR18}
A. Bahrouni and V. D. R\u{a}dulescu,
{\sl On a new fractional Sobolev space and applications to nonlocal variational problems with variable exponent},
Discrete Contin. Dyn. Syst. Ser. S $\mathbf{11}$:3 (2018), 379-389.

\bibitem{T11T}
T. Bartolo, V. Benci and D. Fortunato,
{\sl Abstract critical point theorems and applications to some nonlinearproblems with ``strong'' resonance at infinity},
Nonlinear Anal. $\mathbf{7}$:9 (1983), 981-1012.

\bibitem{3}	
L. Caffarelli,
{\sl Nonlocal equations, drifts and games}, Nonlinear Partial Differential Equations,
 Abel Symposia $\mathbf{7}$ (2012), 37-52.
 
\bibitem{Chung13}
 N.T. Chung,
{\sl Multiple solutions for a $p(x)$-Kirchhoff-type equation with sign-changing nonlinearities},
Complex Var. Elliptic Equ. $\mathbf{58}$:12 (2013), 1637-1646.

\bibitem{CN10}
N.T. Chung and Q.A. Ngo,
{\sl Multiple solutions for a class of quasilinear elliptic equations of $p(x)$-Laplacian type with nonlinear boundary conditions},
Proc. Royal Soc. Edinburgh Sect. A: Math. $\mathbf{140}$:2 (2010), 259-272.

\bibitem{CP}
 F. Colasuonno and P. Pucci,
{\sl Multiplicity of solutions for $p(x)$-polyharmonic Kirchhoff equations},
Nonlinear Anal. $\mathbf{74}$ (2011), 5962-5974.

\bibitem{DM}
G. Devillanova and C.G. Marano, {\sl A free fractional viscous oscillator as a forced standard damped vibration},
Fract. Calc. Appl. Anal. {\bf 19}:2 (2016), 319-356.

\bibitem{DS1}
G. Devillanova and S. Solimini, {\sl Infinitely many positive solutions to some nonsymmetric scalar field equations:
the planar case}, Calc. Var. {\bf 52}:3-4 (2015), 857-898.

\bibitem{DS2}
G. Devillanova and S. Solimini,
 {\sl Some remarks on profile decomposition theorems}, 
 Adv. Nonlinear Stud. {\bf 16}:4 (2016), 795-805.
 
\bibitem{4}	
E. Di Nezza, G. Palatucci and  E. Valdinoci,
{\sl Hitchhiker's guide to the fractional Sobolev spaces}, Bull. Sci. Math. $\mathbf{136}$:5 (2012), 521-573.

\bibitem{5}
L. Diening, P. Harjulehto, P. H\"ast\"o and M. R\r{u}\v{z}i\v{c}ka,
{\sl Lebesgue and Sobolev Spaces with Variable Exponents}, Lecture Notes in Mathematics {\bf 2017},
Springer-Verlag, Heidelberg, 2011.

\bibitem{P3}
A. Fiscella, P. Pucci and B.L. Zhang,
{\sl $p$-fractional Hardy-Schr\"{o}dinger-Kirchhoff systems with critical nonlinearities,}
Adv. Nonlinear Anal. $\mathbf{8}$ (2019), 1111-1131.

\bibitem{H2}
M.K. Hamdani,
{\sl On a nonlocal asymmetric Kirchhoff problems},
Asian-European J. Math.  $\mathbf{13}$:5 (2020), art. 2030001.  

\bibitem{HCR}
M.K. Hamdani, N.T. Chung and D.D. Repov\v{s},
{\sl New class of sixth-order nonhomogeneous $p(x)$-Kirchhoff problems with sign-changing weight functions},
submitted.

\bibitem{Ham}
M.K. Hamdani, A. Harrabi, F. Mtiri and D.D. Repov\v{s},
{\sl Existence and multiplicity results for a new $p(x)$-Kirchhoff problem},
Nonlinear Anal. $\mathbf{190}$ (2020), art. 111598.

\bibitem{HR}
M.K. Hamdani and D.D. Repov\v{s},
{\sl Existence of solutions for systems arising in electromagnetism},
J. Math. Anal. Appl. {\bf  486}:2 (2020), art. 123898.

\bibitem{HHS}
A. Harrabi, M.K. Hamdani, A. Selmi,
{\sl Existence results of the zero mass polyharmonic system},
Complex Var. Elliptic Equ., published online, https://www.tandfonline.com/doi/full/10.1080/17476933.2019.1679794. 

\bibitem{HR16}
Peter H\"{a}st\"{o} and A.M. Ribeiro,
{\sl Characterization of the variable exponent Sobolev norm without derivatives},
Comm. Contemp. Math. $\mathbf{19}$:3 (2017), art. 1650022.

\bibitem{kkk} 
G. Kirchhoff,
Mechanik,  Teubner, Leipzig, 1883.

\bibitem{39}
O. Kov\'{a}\v{c}ik and J. R\'{a}kosn\'{i}k,
{\sl On spaces $L^{p(x)}$ and $W^{1,p(x)}$, }
Czechoslovak Math. J. $\mathbf{41}$:5 (1991), 592-618.

\bibitem{1} J. I. Lee, J. Kim, Y. Kim and J. Lee,
{\sl Multiplicity of weak solutions to non-local elliptic equations involving the fractional $p(x)$-Laplacian},
J. Math. Phys. $\mathbf{61}$:1 (2020), 011505.

\bibitem{LRDZ}
G. Li, V.D. R\u{a}dulescu, D.D. Repov\v{s} and Q. Zhang,
{\sl Nonhomogeneous Dirichlet problems without the Ambrosetti-Rabinowitz condition},
Topol. Methods Nonlinear Anal. $\mathbf{51}$:1 (2018), 55-77.

\bibitem{MBTZ}
X. Mingqi, G. Molica Bisci, G. Tian and B. Zhang,
 {\sl Infinitely many solutions for the stationary Kirchhoff
problems involving the fractional $p$-Laplacian}, 
Nonlinearity {\bf 29} (2016), 357-374.

\bibitem{BRS}
G. Molica Bisci, V. Radulescu and R. Servadei, 
 Variational Methods for Nonlocal Fractional Problems,
 Encyclopedia of Mathematics and its Applications {\bf 162},
Cambridge University Press, Cambridge, 2015.

\bibitem{PRR} 
N.S. Papageorgiou, V.D. R\v{a}dulescu and D.D. Repov\v{s},
 Nonlinear Analysis - Theory and Methods,
  Springer, Cham, 2019.	

\bibitem{P8}
P. Pucci and S. Saldi,
{\sl Critical stationary Kirchhoff equations in $\mathbb{R}^{N}$ involving non-local operators},
Rev. Mat. Iberoam. $\mathbf{32}$ (2016), 1-22.

\bibitem{P9} P. Pucci, M.Q. Xiang and B. Zhang,
{\sl Multiple solutions for nonhomogeneous Schr\"{o}dinger-Kirchhoff type equations involving the fractional $p$-Laplacian in $\mathbb{R}^{N}$},
 Calc. Var. $\mathbf{54}$ (2015), 2785-2806.
 
\bibitem{RD}
V.D. R\u{a}dulescu and D.D. Repov\v{s},
{\sl Partial Differential Equations with Variable Exponents: Variational Methods and Qualitative Analysis},
Taylor \& Francis, CRC Press, Boca Raton, 2015.

\bibitem{R3}
M. Ru\v{z}i\v{c}ka,
\emph{Electro-Rheological Fluids: Modeling and Mathematical Theory}, Lecture Notes in Math. {\bf 1784},
Springer, Berlin, 2000.

\bibitem{T11CC}
X. H. Tang,
 \emph{Infinitely many solutions for semilinear Schr\"{o}dinger equations with sign-changing potential and nonlinearity},
J. Math. Anal. Appl. $\mathbf{401}$ (2013), 407-415.

\bibitem{P12}
M.Q. Xiang,  B. Zhang and  M. Ferrara,
{\sl Existence of solutions for Kirchhoff type problem involving the non-local fractional $p$-Laplacian,}
 J. Math. Anal. Appl. $\mathbf{424}$ (2015), 1021-1041.

\bibitem{P13}
M.Q. Xiang, B. Zhang and V.D. R\u{a}dulescu,
{\sl Multiplicity of solutions for a class of quasilinear Kirchhoff system involving the fractional $p$-Laplacian,}
 Nonlinearity $\mathbf{29}$ (2016), 3186-3205.

\bibitem{ZZ}
C. Zhang and X. Zhang,
{\sl Renormalized solutions for the fractional $p(x)$-Laplacian equation with $L^1$ data},
Nonlinear Anal. $\mathbf{190}$ (2020), art. 111610.

\bibitem{004}
J. Zuo, T. An and M. Li,
{\sl Superlinear Kirchhoff-type problems of the fractional $p$-Laplacian without the (AR) condition},
Bound. Value. Probl. $\mathbf{2018}$ (2018): 180.

\bibitem{ZUO1}
J. Zuo, T. An, L. Yang and X.  Ren,
{\sl The Nehari manifold for a fractional $p$-Kirchhoff system involving sign-changing weight function and concave-convex nonlinearities}, 
J. Funct. Spaces $\mathbf{2019}$ (2019), art. ID 7624373.

\bibitem{005}
J. Zuo, T. An, G. Ye and Z. Qiao,
{\sl Nonhomogeneous fractional $p$-Kirchhoff problems involving a critical nonlinearity},
Electron. J. Qual. Theory Differ. Equ. $\mathbf{41}$ (2019), 1-15.

\end{thebibliography}
\end{document}